\documentclass[production,11pt]{jsg}
\usepackage{amsmath}
\usepackage{amsthm}
\usepackage{amssymb}
\usepackage{amscd}
\usepackage{graphicx}
\usepackage{epsfig}
\input xy
\xyoption{all}
\newtheorem{theorem}{Theorem}[section]

\newtheorem{proposition}[theorem]{Proposition}

\theoremstyle{definition}

\theoremstyle{remark}

\numberwithin{equation}{section}
\allowdisplaybreaks
\begin{document}
%
%
%
\hyphenation{-ho-mo-to-pies}
\hyphenation{ho-mo-to-pies}
\hyphenation{grou-po-ids}
\hyphenation{co-i-so-tro-pic}
\hyphenation{co-i-so-tro-pe}
\hyphenation{con-stra-i-nt}
\hyphenation{sub-ma-ni-fold}
\hyphenation{Mac-ken-zie}
\hyphenation{group-oid}
\newcommand{\bea}{\begin{eqnarray*}}
\newcommand{\eea}{\end{eqnarray*}}
\newcommand{\bdm}{\begin{displaymath}}
\newcommand{\edm}{\end{displaymath}}
\newcommand{\barr}{\begin{array}}
\newcommand{\earr}{\end{array}}
\newcommand{\ben}{\begin{enumerate}}
\newcommand{\een}{\end{enumerate}}
\newcommand{\bde}{\begin{description}}
\newcommand{\ede}{\end{description}}
\newcommand{\nn}{\nonumber}
\newcommand{\nobr}{\nolinebreak}
\newcommand{\noi}{\noindent}
\newcommand{\fn}[1]{\footnote{#1}}
\newcommand{\hbxs}[1]{\hbox{\small{$#1$}}}
\newcommand{\hs}[1]{\hspace{#1cm}}
\newcommand{\vs}[1]{\vspace{#1cm}}
\def\spa{\vskip0.25cm\noindent}
\newcommand{\CC}{\mathbb{C}}
\newcommand{\RR}{\mathbb{R}}
\newcommand{\QQ}{\mathbb{Q}}
\newcommand{\ZZ}{\mathbb{Z}}
\newcommand{\NN}{\mathbb{N}}
\newcommand{\KK}{\mathbb{K}}
\newcommand{\HH}{\mathbb{H}}
\newcommand{\FF}{\mathbb{F}}
\newcommand{\Z}[1]{\mathbb{Z}_{#1}}
\newcommand{\ds}{\mathbb{S}}
\newcommand{\BB}{\mathbb{B}}
\newcommand{\DD}{\mathbb{D}}
\newcommand{\calA}{\mathcal{A}}
\newcommand{\calB}{\mathcal{B}}
\newcommand{\calC}{\mathcal{C}}
\newcommand{\calD}{\mathcal{D}}
\newcommand{\calF}{\mathcal{F}}
\newcommand{\calH}{\mathcal{H}}
\newcommand{\calV}{\mathcal{V}}
\newcommand{\calS}{\mathcal{S}}
\newcommand{\calM}{\mathcal{M}}
\newcommand{\calN}{\mathcal{N}}
\newcommand{\calK}{\mathcal{K}}
\newcommand{\calG}{\mathcal{G}}
\newcommand{\calI}{\mathcal{I}}
\newcommand{\calJ}{\mathcal{J}}
\newcommand{\calL}{\mathcal{L}}
\newcommand{\call}{\ell}
\newcommand{\calO}{\mathcal{O}}
\newcommand{\calR}{\mathcal{R}}
\newcommand{\tcalG}{\Tilde{\mathcal{G}}}
\newcommand{\ttcalG}{\Tilde{\tcalG}}
\newcommand{\ucalG}{\underline{\mathcal{G}}}
\newcommand{\bucalG}{\boldsymbol{\ucalG}}
\newcommand{\calE}{\mathcal{E}}
\newcommand{\calP}{\mathcal{P}}
\newcommand{\calW}{\mathcal{W}}
\let\sf=\mathsf
\let\tsf=\textsf
\newcommand{\etsf}[1]{\tsf{\emph{#1}}}
\newcommand{\sfA}{\mathsf{A}}
\newcommand{\sfB}{\mathsf{B}}
\newcommand{\sfD}{\mathsf{D}}
\newcommand{\sfH}{\mathsf{H}}
\newcommand{\sfS}{\mathsf{S}}
\newcommand{\sfC}{\mathsf{C}}
\newcommand{\sfM}{\mathsf{M}}
\newcommand{\sfK}{\mathsf{K}}
\newcommand{\sfG}{\mathsf{G}}
\newcommand{\sfI}{\mathsf{I}}
\newcommand{\sfL}{\mathsf{L}}
\newcommand{\sfO}{\mathsf{O}}
\newcommand{\tsfG}{\Tilde{\mathsf{G}}}
\newcommand{\ttsfG}{\Tilde{\tsfG}}
\newcommand{\usfG}{\underline{\mathsf{G}}}
\newcommand{\busfG}{\boldsymbol{\usfG}}
\newcommand{\sfE}{\mathsf{E}}
\newcommand{\sfP}{\mathsf{P}}
\let\rm=\mathrm
\let\trm=\textrm
\newcommand{\rmA}{\mathrm{A}}
\newcommand{\rmB}{\mathrm{B}}
\newcommand{\rmH}{\mathrm{H}}
\newcommand{\rmS}{\mathrm{S}}
\newcommand{\rmC}{\mathrm{C}}
\newcommand{\rmM}{\mathrm{M}}
\newcommand{\rmK}{\mathrm{K}}
\newcommand{\rmG}{\mathrm{G}}
\newcommand{\rmI}{\mathrm{I}}
\newcommand{\rmL}{\mathrm{L}}
\newcommand{\rmO}{\mathrm{O}}
\newcommand{\trmG}{\Tilde{\mathrm{G}}}
\newcommand{\ttrmG}{\Tilde{\trmG}}
\newcommand{\urmG}{\underline{\mathrm{G}}}
\newcommand{\burmG}{\boldsymbol{\urmG}}
\newcommand{\rmE}{\mathrm{E}}
\newcommand{\rmP}{\mathrm{P}}
\newcommand{\bfA}{\mathbf{A}}
\newcommand{\bfB}{\mathbf{B}}
\newcommand{\bfH}{\mathbf{H}}
\newcommand{\bfS}{\mathbf{S}}
\newcommand{\bfC}{\mathbf{C}}
\newcommand{\bfM}{\mathbf{M}}
\newcommand{\bfN}{\mathbf{N}}
\newcommand{\bfK}{\mathbf{K}}
\newcommand{\bfG}{\mathbf{G}}
\newcommand{\bfI}{\mathbf{I}}
\newcommand{\bfL}{\mathbf{L}}
\newcommand{\bfO}{\mathbf{O}}
\newcommand{\tbfG}{\Tilde{\mathbf{G}}}
\newcommand{\ttbfG}{\Tilde{\tbfG}}
\newcommand{\ubfG}{\underline{\mathbf{G}}}
\newcommand{\bubfG}{\boldsymbol{\ubfG}}
\newcommand{\bfE}{\mathbf{E}}
\newcommand{\bfP}{\mathbf{P}}
\newcommand{\bfj}{\mathbf{j}}
\newcommand{\bfs}{\mathbf{s}}
\newcommand{\bft}{\mathbf{t}}
\def\id{{\sf{id}}}
\def\ker{{\sf{ker}\,}}
\def\Ker{{\pmb{\mathcal{K}\emph{er}}}\,}
\def\im{{\sf{Im}}}
\def\pr{{\mathrm{pr}}}
\def\Pr{{\mathrm{Pr}}}
\def\rank{{\sf{rank}\,}}
\def\rk{{\sf{rk}\,}}
\def\mon{{\sf{Mon}}}
\def\fii{\varphi}
\def\eps{\varepsilon}
\let\Hat=\widehat
\def\comp{
\circ
}
\newcommand{\smalcirc}{\mbox{\tiny{$\circ $}}} 
\def\inverse{^{\hbox{-\tiny{$1$}}}}
\def\wt{\widetilde}
\def\mod{/\hbox{\footnotesize{$\thicksim$}}}
\newcommand{\tover}[1]{\times_{{#1}}\!}
\newcommand{\twice}[1]{{#1}^{\times 2}}
\newcommand{\threetimes}[1]{{#1}^{\times 3}}
\def\acts{\ltimes
}
\def\lrtimes{{\ltimes\hs{-0.35}\rtimes}} 
\newcommand{\ra}{\rangle}
\newcommand{\la}{\langle}  
\def\inc{\hookrightarrow}
\def\sur{\twoheadrightarrow}
\let\rinv=\overrightarrow
\let\linv=\overleftarrow
\let\lra=\longrightarrow
\newcommand{\rrinv}[1]{\overset\Rightarrow{#1}}
\newcommand{\llinv}[1]{\overset\Rightarrow{#1}}

\def\pa{\partial}
\def\dd{\mathrm{d}}
\newcommand{\pad}[2]{\frac{\pa{#1}}{\pa{#2}}}
\newcommand{\tod}[2]{\frac{\dd{#1}}{\dd{#2}}}
\newcommand{\vad}[2]{\frac{\delta{#1}}{\delta{#2} }}
\newcommand{\ddt}[1]{\frac{\dd{#1}}{\dd{t}}}
\def\cif{\calC^\infty}
\def\frax{\mathfrak{X}}
\def\ham{{\sf{Ham}}}
\def\osharp{\omega^\sharp}
\def\oosharp{\Omega^\sharp}
\def\pish{\pi^\sharp}
\def\ppish{\Pi^\sharp}
\def\ohm{\omega}
\def\Ohm{\Omega}
\def\frag{\mathfrak{g}}
\def\frah{\mathfrak{h}}
\def\frack{\mathfrak{k}}
\def\frai{\mathfrak{i}}
\def\frak{\mathfrak{c}}
\def\gpd{\,\lower1.5pt\hbox{$\rightarrow$}\hskip-3.5mm\raise2.5pt\hbox{$\rightarrow$}\,}       
\def\dpg{\,\lower1pt\hbox{$\leftarrow$}\hskip-3.75mm\raise3pt\hbox{$\leftarrow$}\,}       
\newcommand{\backl}{\mathbin{\vrule width1.5ex height.4pt\vrule height1.5ex}}
\newcommand{\per}{\backl}     
\newcommand{\poidd }[2]{#1\gpd #2}
\newcommand{\poiddd }[3]{ (#1\gpd #2, \alpha_{#3}, \beta_{#3})}    
\newcommand{\gammar}[1]{\rinv{\Gamma}{(#1)}}
\newcommand{\gammal}[1]{\linv{\Gamma}{(#1)}}
\newcommand{\gammam}[1]{\widetilde{\Gamma}{(#1)}}
\def\rgamma{\rinv{\Gamma}{(\calG,\Ohm)}}
\def\lgamma{\linv{\Gamma}{(\calG,\Ohm)}}
\def\mgamma{\widehat{\Gamma}{(\calG,\Ohm)}}
\def\pgamma{\Gamma^{\downarrow}{(\calG,\Ohm)}}
\def\qgamma{\Gamma(M/\calG,A/\Ohm)}
\def\mmgamma{\Hat{\Gamma}^\bullet(M,\calG)}
\def\pgammaa{\Gamma^{\downarrow}}
\def\pfrax{\frax^{\downarrow}}
\newcommand{\brak}[2]{[{\,#1\,,\,#2\,}]}
\newcommand{\bracts}[2]{[{\,#1\,\overset{\acts},\,#2\,}]}
\newcommand{\brastar}[2]{[{#1
{\overset\star{,}}\,
#2}]}
\newcommand{\brast}[2]{[{#1
{\overset{*}{,}}\,
#2}]}
\newcommand{\poib}[2]{\{{\,#1\,,\,#2\,}\}}
\newcommand{\pair}[2]{\langle{\,#1\,,\,#2\,}\rangle}

\def\Dots{,\,\dots\,,}


\newcommand{\pmbgpdm}{$\poidd{\pmb\calG}{\bf M}\:\:$}
\newcommand{\pmbact}{$\poidd{\pmb\calG\acts\bfN}{\bf N}\:\:$}
\def\pmbs{\pmb\sigma}
\def\pmbd{\pmb\delta}
\def\pmbk{\pmb\chi}
\def\pmbi{\pmb\iota}
\def\pmbe{\pmb\eps}
\def\pmbg{\hbox{$\pmb\calG$}}

\let\Bar=\overline
\let\ol=\overline
\let\ul=\underline

\let\bgn=\begin

\newcommand{\temise}[1]{{
\begin{itemize}
{#1}
\end{itemize}
}}
\def\it{\item[$i$)]}
\def\iit{\item[$ii$)]}
\def\iiit{\item[$iii$)]}
\def\ivt{\item[$iv$)]}
\def\vt{\item[$v$)]}
\def\vit{\item[$vi$)]}
\def\viit{\item[$vii$)]}
\def\viiit{\item[$viii$)]}
\def\punto{\item[$\cdot$]}

\def\ut{\item[1)]}
\def\dt{\item[2)]}
\def\tt{\item[3)]}
\def\qt{\item[4)]}
\def\ct{\item[5)]}
\def\st{\item[6)]}

\newcommand{\be}{\begin{eqnarray*}}
\newcommand{\ee}{\end{eqnarray*}}  

\newcommand{\Exasss}[4]{{\Exa{{#1}\temise{
\item[1.] {#2}
\item[3.] {#3}
\item[4.] {#4}
}}}}

\def\Lg{\mathcal{LG}}
\def\La{\mathcal{LA}}
\def\tlg{$\Lg$}
\def\tla{$\La$}
\def\tvb{$\mathcal{VB}$}

\def\daga{^{\hbox{\tiny{$\pmb{++}$}}}\!}
\def\ppv{{\iota_{\hbox{\tiny{$\pa_V^+$}}}}^{\hs{-0.315}*}\hs{0.2}}
\def\pmv{{\iota_{\hbox{\tiny{$\pa_V^-$}}}}^{\hs{-0.315}*}\hs{0.2}}
\def\pph{{\iota_{\hbox{\tiny{$\pa_H^+$}}}}^{\hs{-0.315}*}\hs{0.2}}
\def\pmh{{\iota_{\hbox{\tiny{$\pa_H^-$}}}}^{\hs{-0.315}*}\hs{0.2}}
\def\ppmv{{\iota_{\hbox{\tiny{$\pa_V^\pm$}}}}^{\hs{-0.315}*}\hs{0.2}}
\def\ppmh{{\iota_{\hbox{\tiny{$\pa_H^\pm$}}}}^{\hs{-0.315}*}\hs{0.2}}

\def\ttar{\mathrm{\hbox{\small{$\mathrm{t}$}}}}
\def\tar{\mathrm{t}}
\def\ssor{\mathrm{\hbox{\small{$\mathrm{s}$}}}}
\def\sor{\mathrm{s}}
\def\vup{\hbox{\tiny{${}_{V}$}}}
\def\hup{\hbox{\tiny{${}_{H}$}}}
\def\vuppp{\hbox{\tiny{${}^{3V}$}}}

\def\ttv{\,{\ttar}\vup}
\def\tsv{\,{\ttar}_{v}}
\def\stv{\,{\ssor}\vup}
\def\ssv{\,{\ssor}_{v}}

\def\sc{{\sor}^c}
\def\tc{{\tar}^c}

\def\tth{\,{\ttar}\hup}
\def\tsh{\,{\ttar}_{h}}
\def\sth{\,{\ssor}\hup}
\def\ssh{\,{\ssor}_{h}}

\def\nttv{{\mathrm{t}}\vup}
\def\ntsv{{\mathrm{t}}^{v}}
\def\nstv{{\mathrm{s}}\vup}
\def\nssv{{\mathrm{s}}^{v}}

\def\ntth{{\mathrm{t}}\hup}
\def\ntsh{{\mathrm{t}}^{h}}
\def\nsth{{\mathrm{s}}\hup}
\def\nssh{{\mathrm{s}}^{h}}

\def\etv{{\eps}\vup}
\def\esv{{\eps}^{v}}
\def\etv{{\eps}\vup}
\def\esv{{\eps}^{v}}
\def\ec{{\eps}^{c}}

\def\eth{{\eps}\hup}
\def\esh{{\eps}^{h}}
\def\eth{{\eps}\hup}
\def\esh{{\eps}^{h}}

\def\itv{{\iota}\vup}
\def\isv{{\iota}^{v}}
\def\itv{{\iota}\vup}
\def\isv{{\iota}^{v}}

\def\ith{{\iota}\hup}
\def\ish{{\iota}^{h}}
\def\ith{{\iota}\hup}
\def\ish{{\iota}^{h}}

\def\mtv{{\mu}\vup}
\def\msv{{\mu}^{v}}
\def\mtv{{\mu}\vup}
\def\msv{{\mu}^{v}}

\def\mth{{\mu}\hup}
\def\msh{{\mu}^{h}}
\def\mth{{\mu}\hup}
\def\msh{{\mu}^{h}}

\def\dtv{{\delta}\vup}
\def\dsv{{\delta}^{v}}
\def\dtv{{\delta}\vup}
\def\dsv{{\delta}^{v}}

\def\dth{{\delta}\hup}
\def\dsh{{\delta}^{h}}
\def\dth{{\delta}\hup}
\def\dsh{{\delta}^{h}}

\def\epsh{\hat{\eps}}
\def\sh{\hat{\sor}}
\def\th{\hat{\tar}}
\def\iotah{\hat{\iota}}
\def\muh{\hat{\mu}}
\def\divh{\hat{\delta}}
\def\rhoh{\hat{\rho}}
\def\chih{\hat{\chi}}
\def\jh{\hat{j\:}}
\def\sigmah{\hat{\sigma}}
\def\msigma{\pmb{\sigma}}

\def\hdot{\,\hat{\cdot}\,}

\def\epsha{\hat{\eps}_{\acts}}
\def\sha{\hat{\sor}_{\acts}}
\def\tha{\hat{\tar}_{\acts}}
\def\iotaha{\hat{\iota}_{\acts}}
\def\muha{\hat{\mu}_{\acts}}
\def\divha{\hat{\delta}_{\acts}}
\def\chiha{\hat{\chi}_{\acts}}

\def\epsa{{\eps}_{\acts}}
\def\sa{{\sor}_{\acts}}
\def\ta{{\tar}_{\acts}}
\def\iotaa{{\iota}_{\acts}}
\def\mua{{\mu}_{\acts}}
\def\diva{{\delta}_{\acts}}
\def\pra{\mathrm{pr}_{\acts}}
\def\Pra{\mathrm{Pr}_{\acts}}

\def\epst{\tilde{\eps}}
\def\st{\tilde{s}}
\def\tt{\tilde{t}}
\def\iotat{\tilde{\iota}}
\def\mut{\tilde{\mu}}

\newcommand{\fib}[2]{\,{}_{#1}\hs{-0.05}\times\hs{-0.05}{}_{#2}\,}
\newcommand{\Lie}[1]{\mathsf{Lie}[#1]}
\newcommand{\lie}[2]{\mathsf{Lie}_{#1}[#2]}

\def\lagpd{\mathcal{LA}\hbox{-groupoid}}
\def\gpdm{\calG\gpd M}
\def\gpdmc{\Tilde{\calG}\gpd M}
\def\pgpd{(\calG,\Pi^G)\gpd M}
\def\pgpdb{(\ol{\calP},{\ol\Pi})\gpd M}
\def\poidm{\calP\gpd M}
\def\poidmb{\bar{\calP}\gpd M}
\def\bagd{(A,A^*)\rightarrow M}

\def\tsp{T^*\calP}
\newcommand{\gr}[1]{\sf\Gamma({#1})}
\newcommand{\ki}[2]{\mathfrak{X}^{#1}{(#2)}}
\def\jm{{\hbox{\j}}}
\def\jmod{\calJ\inverse(e_\star)/G}
\def\jj{\calJ\inverse(e_\star)}
\def\sg{\hbox{\tiny{$G$}}}
\def\dim{\sf{dim}\,}



\markboth{Luca Stefanini}{
Integrability and reduction of Poisson group actions
}
$ $
\bigskip

\bigskip

\centerline{{\large{\textsc{
Integrability and reduction of Poisson group actions
}}
}\footnote{\emph{2000 Mathematics Subject Classification}:
Primary 53D20; Secondary 58H05, 18D05.}}

\bigskip
\centerline{{\large by Luca Stefanini}}
\bigskip
{\small
\emph{
\bgn{center}{
Institut f\"ur Mathematik, Universit\"at Z\"urich-Irchel,\\
Winterthurerstrasse 190, CH-8057 Z\"urich, Switzerland}
\end{center}}
\centerline{\tsf{lucaste@math.unizh.ch} and \tsf{cucanini@gmail.com}}}
\medskip

\centerline{30.10.2007}

\medskip
\subsection*{Abstract}
In this paper we study Poisson actions of complete  Poisson groups,
without any connectivity assumption or requiring the existence of a
momentum map. For any complete Poisson group $G$  with dual $G^\star$
we obtain a suitably connected integrating symplectic double groupoid
$\calS$. As a consequence, the cotangent lift of a Poisson action on
an integrable Poisson manifold $P$ can be integrated to a Poisson
action of the symplectic groupoid $\poidd{\calS}{G^\star}$ on the
symplectic groupoid for $P$. Finally, we show that the quotient
Poisson manifold $P/G$ is also integrable, giving an explicit
construction of a symplectic groupoid for it, by a reduction
procedure on an associated morphism of double Lie groupoids.  
\tableofcontents
\section*{Introduction}
Symplectic groupoids were introduced independently by Karas\" ev
\cite{ka}, Weinstein \cite{ws} and Zakrzewski \cite{za} in the
attempt of providing a geometric setting for the quantization of
Poisson manifolds. The symplectic foliation of an integrable Poisson
bivector is described by the orbits of an integrating symplectic
groupoid, therefore symplectic groupoids also represent a fundamental
tool in the study of the classical geometry of Poisson manifolds.
\\
In recent years, the integration problem for Poisson manifolds to
symplectic groupoids has been completely solved. Mackenzie and Xu
showed in \cite{mx2} that a Poisson manifold $P$ is integrable to a
symplectic groupoid if and only if the associated Lie algebroid on
$T^*P$ is integrable. Cattaneo and Felder produced in \cite{cafe} a
topological model for the symplectic groupoid of \emph{any} Poisson
manifold, yielding a smooth symplectic groupoid in the integrable
case, by the symplectic reduction of the Poisson sigma model. Finally
Crainic and Fernandes obtained in \cite{crfs1,crfs2} necessary and
sufficient conditions for a Poisson manifold to be integrable, i.e.
to endow the topological groupoid of Cattaneo and Felder with a
smooth structure.
\\
On the one hand Poisson manifolds behave well under reduction for
actions as general as those of Poisson groupoids: a Poisson bivector
always descends to the quotient by a free and proper compatible
action \cite{wsc}. On the other hand, Poisson manifolds are not
always integrable to symplectic groupoids; the natural question is
thus
\bgn{quote}
{\em
Are quotients of integrable Poisson manifolds also integrable?
}
\end{quote}
A positive answer was recently given by Fernades-Ortega-Ratiu
\cite{for} in the case of Lie group actions and Lu \cite{luph} gave a
construction  of symplectic groupoids for certain Poisson homogeneous
spaces. In the case of Poisson group actions with a complete momentum
mapping, Xu described in \cite{xu} a reduction procedure on a lifted
moment map; when the reduced space of the latter is smooth, it is a
symplectic groupoid for the quotient Poisson structure.
\\
In this paper we give a positive answer to the above question in the
case of free and proper Poisson actions of complete Poisson groups.
Under these assumptions on a Poisson $G$-space $P$, with symplectic
groupoid $\Lambda$, the integrability of the left  dressing action
$\frag^*\to\frax(G)$ allows to lift the action of $G$ on $P$ to the 
action of  suitably connected symplectic double groupoid on a
canonical morphism of Poisson groupoids $\calJ:\Lambda\to G^\star$ to
the dual Poisson group. Moreover $\calJ$ extends to a ``moment
morphism'' of double Lie groupoids, whose kernel (double Lie
groupoid) provides the  data for a reduction {\em \`a la}
Mikami-Weinstein \cite{MiWs88} to produce an integration of the
quotient Poisson bivector on $P/G$.
\smallskip
Our main result may be rephrased as the following
\bgn{thm}[\bf \ref{redux}] Let $G$ be a complete Poisson group and
$P$ a Poisson $G$-space. If $G$ acts freely and properly and $P$ is
integrable, the quotient Poisson manifold $P/G$ is also integrable;
in particular, the quotient $\jmod$ is smooth and the unique
symplectic form  $\ul{\ohm}$  on $\jmod$  such that
$\pr^*\ul{\ohm}=\iota^*\ohm$ makes it a symplectic groupoid for ther
quotient Poisson bivector on $P/G$.
\end{thm}
\noindent
On the way to our main result we also show that any complete Poisson
group (not necessarily 0- or 1-connected) is integrable to a
symplectic double groupoid (theorem \ref{main}), extending a result
by Lu and Weinstein \cite{lw}. Moreover, we also produce a class of
examples for which our functorial approach  to the integrability of
Lie algebroids \cite{st} is effective; namely, all action
\tla-groupoids associated with Poisson actions of complete Poisson
groups on integrable Poisson manifolds are integrable by suitable
action double Lie groupoids (theorem \ref{mainapp}). 
\smallskip
Mackenzie's approach to Poisson reduction \cite{mu} has been a major
source of inspiration for this paper; in fact, we develop here an
integrated version of the reduction procedure described in \cite{mu}.
\smallskip
The reduction of morphisms of Poisson groupoids by compatible actions
of symplectic double groupoids can be performed along the lines of
this paper. This approach is however ineffective to obtain symplectic
groupoids for quotients of Poisson group\emph{oid} actions. We shall
discuss in a forthcoming paper the integrability of this type of
quotient Poisson manifolds.
\smallskip
This paper is organized as follows.
The first two sections are introductory.
In \S\ref{dpg} we present the main definitions and review the notion
of duality for Poisson groupoids in relation with integrability
issues%
\footnote{ 
In this paper, apart from \S\ref{dpg}, we shall only consider dual
Poisson groups.  We chose to define duality in the context of Poisson
groupoids, since it allows us to introduce quickly and more naturally
\tla-groupoids and double Lie groupoids, which we use throughout the
rest of the paper.
}.
In \S\ref{llah} we fix our conventions for Lie algebroid homotopies
and recall the lifting conditions of \cite{st} for the integrability
of fibred products of Lie algebroids.
In \S\ref{dasdg} we discuss the integrability of complete Poisson
groups to symplectic double groupoids.
In  \S\ref{ipgs} we lift the action of a complete Poisson group to a 
compatible action of the associated symplectic double groupoid.
In \S\ref{iqps} we perform the reduction of the ``moment
morphism''.   We finally comment on Xu's reduction of \cite{xu} and
Fernandes-Ortega-Ratiu's  integration of \cite{for}, which we obtain
as special cases. In particular our construction shows that Xu's
quotient is always smooth; we conclude by remarking that the
condition of \cite{for} for integration to commute with reduction
extends to our setting.
\smallskip
\noindent
We fix at the end of this paper some of our conventions and special
notations.
\section{Duality for Poisson groupoids}\label{dpg}
A \tsf{Poisson groupoid} is a Lie groupoid $\poidd{\calG}{M}$ endowed
with a Poisson bivector $\Pi$ making  the graph $\gr{\mu}$ of the
partial multiplication  $\mu:\calG^{(2)}\to \calG$ coisotropic in
$\calG\times\calG\times\ol{\calG}$. When $\Pi$ is non degenerate,
hence $\Pi^\sharp =\ohm^{\sharp-1}$ for some symplectic form $\ohm$,
the compatibility condition is equivalent to $\gr{\mu}$ being
Lagrangian and $\calG$ is called a \tsf{symplectic groupoid}.  A
Poisson groupoid over the one point  manifold is a \tsf{Poisson
group} and its multiplication is a Poisson map. Any Lie group(oid) is
a Poisson group(oid) for the zero Poisson bivector.
\\ 
The main properties of Poisson groupoids were proved in \cite{wsc};
in particular the unit section $M\to\calG$ is a coisotropic embedding
(a Lagrangian embedding in the symplectic case) and there exists a
unique Poisson structure on $M$ making the source map Poisson. We
shall say that the special Poisson structure on the base of a 
Poisson groupoid is the \tsf{induced Poisson structure}.  Consider
now any Poisson manifold $(P,\pi)$; there exists a unique Lie
algebroid bracket on $\Ohm^1(P)$ extending the Lie bracket
$$
\qquad
\brak{\dd f_+}{\dd f_-}:=\dd \poib{f_+}{f_-}\qquad,\qquad f_\pm\in\cif(P)\qquad,
$$
on the space of exact 1-forms according to the Leibniz rule, that is
compatible with the Poisson anchor $\pish:T^*P\to TP$. Whenever such
a Lie algebroid is integrable to a Lie groupoid 
$\poidd{\Lambda}{P}$, one can show \cite{mx2} that $\Lambda$,
provided it is source 1-connected, carries a unique symplectic form
making it a symplectic groupoid and inducing $\pi$ on the base
manifold $P$. In that case $P$ is called an \tsf{integrable Poisson
manifold} and $\Lambda$ is called \emph{the} symplectic groupoid  of
$P$.
\\
Let us describe a typical example of an integrable Poisson manifold.
For any Lie group $G$ with Lie algebra $\frag$, $T^*G$ can be endowed
with a Lie groupoid structure over $\frag^*$ making it a symplectic
groupoid for the linear Poisson bivector on $\frag^*$ dual to the Lie
bracket of $\frag$ %
 (the construction applies, mutatis mutandis, to any Lie groupoid
$\calG$ with Lie algebroid $A$ to produce a symplectic groupoid
$\poidd{T^*\calG}{A^*}$, see \cite{cdw} for details). %
 Source $\sh$ and target $\th$ are the evaluations of the left and
right vector bundle trivializations $T^*G\to \frag^*\times G$, the
inversion is the anti-transpose $\iotah=-\dd\,\iota^{t}$ of the
tangent group inversion and the unit  section $\epsh$ is the
identification of $\frag^*$ with the conormal bundle $N^*{e}$ of the
identity element $e\subset G$. One can check that the conormal 
bundle $N^*\gr{\mu}\subset T^*G^{\hbox{\tiny{$\times 3$}}}$ of the
group multiplication is the graph of a map, therefore setting
$\gr{\muh}:=(\id_{T^*G}\times \id_{T^*G}\times
-\id_{T^*G})N^*\gr{\mu}$, yields a groupoid multiplication
$\muh:(T^*G)^{(2)}\to T^*G$, whose graph is Lagrangian in 
$T^*G\times T^*G\times \ol{T^*G}$, by construction, therefore making
making $\poidd{T^*G}{\frag^*}$ symplectic groupoid. The  Poisson
bracket on $\frag^*$ is induced by the anti-canonical symplectic
form.
Let now $(G,\Pi)$ be a Poisson group. The linearization of $\Pi$ at
the identity endows $\frag^*$ with a compatible Lie bracket
$\brast{\:\,}{\,}$,
$$
\qquad
\pair{\brast{\xi_+}{\xi_-}}{x}:=(\pounds_{\linv{x}}\Pi)_e(\xi_+,\xi_-)
\qquad,\qquad \xi_\pm\in\frag^*\quad,\quad x\in\frag
\qquad,
$$
for the left invariant vector field  $\linv{x}$ of $x\in\frag$.  Such
a pair $(\frag,\frag^*)$ is called a \tsf{Lie bialgebra}.  We shall
remark that the direct sum $\frag\oplus\frag^*$ has a natural Lie
algebra structure,  the \tsf{Drinfel'd double} $\mathfrak{d}$,
obtained as a double twist of the canonical bracket on
$\frag\oplus\frag^*$,  making $\frag$ and $\frag^*$ Lie subalgebras
and the canonical Riemannian structure  ad-invariant%
\footnote{ 
For an account of the basic facts on Poisson groups and Lie
bialgebras  see \cite{cp}.
}.
\\ 
More generally Poisson groupoids differentiate to Lie bialgebroids: a
\tsf{Lie bialgebroid $(A,A^*)$} \cite{mx1}  is a pair of  Lie
algebroids in duality such that the Lie algebroid differential
$\Gamma(\wedge^\bullet A)\to \Gamma(\wedge^{\bullet +1} A)$ induced
by $A^*$ is a derivation of the graded Lie bracket on 
$\Gamma(\wedge^\bullet A)$ induced by $A$ \cite{ks}.  The source
1-connected integration of either Lie algebroid of a Lie bialgebroid,
when it exists, can be endowed with a unique Poisson  structure
making it a Poisson groupoid and inducing the given Lie bialgebroid
\cite{mx2}.
$$
\bgn{array}{ccccccc}
\xy
*+{}="0",    <-0.6cm,0.6cm>
*+{\calD}="1", <0.6cm,0.6cm>
*+{\calV}="2", <-0.6cm,-0.6cm>
*+{\calH}="3", <0.6cm,-0.6cm>
*+{M}="4",
\ar  @ <0.07cm>   @{->} "1";"2"^{} 
\ar  @ <-0.07cm>  @{->} "1";"2"_{}  
\ar  @ <-0.07cm>  @{->} "1";"3"_{}
\ar  @ <0.07cm>   @{->} "1";"3"^{}
\ar  @ <0.07cm>   @{->} "2";"4"^{}
\ar  @ <-0.07cm>  @{->} "2";"4"_{}
\ar  @ <-0.07cm>  @{->} "3";"4"_{}
\ar  @ <0.07cm>   @{->} "3";"4"^{}
\endxy
&\quad\:\:&
\xy
*+{}="0",    <-0.6cm,0.6cm>
*+{\Omega}="1", <0.6cm,0.6cm>
*+{A}="2", <-0.6cm,-0.6cm>
*+{\calG}="3", <0.6cm,-0.6cm>
*+{M}="4",
\ar  @ <-0.07cm>   @{->} "1";"2"^{} 
\ar  @ <0.07cm>    @{->} "1";"2"^{}  
\ar  		   @{->} "1";"3"^{}
\ar                @{->} "2";"4"_{}
\ar  @ <0.07cm>    @{->} "3";"4"^{}
\ar  @ <-0.07cm>   @{->} "3";"4"_{}
\endxy
&\quad\:\:&
\xy
*+{}="0",    <-0.6cm,0.6cm>
*+{T^*G}="1", <0.6cm,0.6cm>
*+{\frag^*}="2", <-0.6cm,-0.6cm>
*+{G}="3", <0.6cm,-0.6cm>
*+{\bullet}="4",
\ar  @ <-0.07cm>   @{->} "1";"2"^{} 
\ar  @ <0.07cm>    @{->} "1";"2"^{}  
\ar  		   @{->} "1";"3"^{}
\ar                @{->} "2";"4"_{}
\ar  @ <0.07cm>    @{->} "3";"4"^{}
\ar  @ <-0.07cm>   @{->} "3";"4"_{}
\endxy
&\quad\:\:&
\xy
*+{}="0",    <-0.6cm,0.6cm>
*+{T^*\calG}="1", <0.6cm,0.6cm>
*+{A^*}="2", <-0.6cm,-0.6cm>
*+{\calG}="3", <0.6cm,-0.6cm>
*+{M}="4",
\ar  @ <-0.07cm>   @{->} "1";"2"^{} 
\ar  @ <0.07cm>    @{->} "1";"2"^{}  
\ar  		   @{->} "1";"3"^{}
\ar                @{->} "2";"4"_{}
\ar  @ <0.07cm>    @{->} "3";"4"^{}
\ar  @ <-0.07cm>   @{->} "3";"4"_{}
\endxy
\\
Figure\: 1
&\quad\:\:&
Figure\:2
&\quad\:\:&
Figure\:3
&\quad\:\:&
{Figure\: 4}
\end{array}
$$
On the other hand Poisson groupoids, which have integrable Poisson
structure,  integrate to symplectic double groupoids, at least in
some favorable cases. A  \tsf{double Lie groupoid} \cite{mdlg} (fig.
1) is a groupoid object in the category of groupoids, such that all
sides are Lie groupoids and the horizontal structural maps are
morphisms of Lie groupoids over the side horizontal maps for the
vertical Lie groupoid structures. In order to make the domains of the
top multiplications Lie groupoids over the domains of the side
multiplications,  the double source map
$
\mathbb{S}:\calD\to\calH\fib{\ssh}{\ssv}\calV
$,$
d\mapsto (\nstv(d), \nsth(d))
$ 
is required to be  submersion%
\footnote{ 
In \cite{mdlg} the double source map is also required to be
surjective. We drop the condition, since it plays a crucial role only
in the theory of transitive double Lie groupoids.
}.
A \tsf{symplectic double groupoid} is a double Lie groupoid endowed
with a symplectic form making both the top groupoids
$\poidd{\calD}{\calV}$ and $\poidd{\calD}{\calH}$  symplectic
groupoids. It follows \cite{msdg} that the induced Poisson structures
on the  side groupoids of a symplectic double groupoid $\calD$ as
above are  Poisson groupoids in duality, in the sense that the
corresponding Lie bialgebroids are isomorphic up to a canonical flip
operation: if $(A,A^*)$ is the Lie bialgebroid of $\calH$, the Lie
bialgebroid of $\calV$ is  canonically isomorphic to $(A^*, -A)$,
where $-A$ is the Lie algebroid obtained by changing the signs of
bracket and anchor of $A$.
For any Poisson groupoid $\poidd{(\calG,\Pi)}{M}$ with integrable
Poisson structure it is easy to see that $A^*$ is an integrable Lie
algebroid, thus the source 1-connected integration
$\poidd{\calG^\star}{M}$ carries a unique Poisson structure
$\Pi^\star$ making it a Poisson groupoid with Lie bialgebroid $(A^*,
-A)$.
Following the conventions of \cite{msdg,st}, we shall say that
$\poidd{(\calG^\star,\Pi^\star)}{M}$ is the (\tsf{weak}) \tsf{dual}
of the given  Poisson groupoid%
\footnote{ 
Declaring the dual groupoid of an integrable Poisson groupoid to be
that with Lie bialgebroid $(A, A^*)$ is also frequent in the
literature.
}.
A \tsf{double} of a Poisson groupoid is a symplectic double groupoid
with the given Poisson groupoid and the unique dual as side Poisson
groupoids: such a double provides a stronger duality relation between
its side groupoids.%
\\
A double Lie groupoid is a groupoid object in the category of Lie
groupoids and applying the Lie functor to the its vertical groupoids
yields an \tsf{\tla-groupoid} \cite{mdlg}, i.e  a diagram such as
that in figure 2, where both horizontal sides are Lie groupoids  and
both vertical sides are Lie algebroids, making the top groupoid
structural maps  morphisms of Lie algebroids over the side ones;
analogously to the case of double Lie groupoids, the double source
map 
$
\$:\Ohm\rightarrow \calG\fib{\sor}{\pr}A
$,
$
\ohm\mapsto (\Pr(\ohm),\sh(\ohm))
$
is required to be surjective, equivalently the top source map to be
fibrewise  surjective, to make the domain of the top multiplication a
Lie algebroid over the domain of the side multiplication. If $G$ is a
Poisson group, the diagram in figure 3 is a typical example of an 
\tla-groupoid. Note that the source regularity condition is fulfilled
since the top source map is fibrewise an isomorphism of vector
spaces.%
\\
Following a functorial approach, one can see \cite{st} that
\tla-groupoids with integrable top Lie algebroid do integrate to
double Lie groupoids, provided some lifting conditions on the top
source map of the \tla-groupoid are met. The result specializes to
the \tla-groupoid (fig. 4) of a Poisson groupoid and one can further
show  that the integrating double Lie groupoid is indeed a double of
the original Poisson groupoid (see \cite{msdg} for the definition of
the \tla-groupoid structure). We shall check in \S\ref{dasdg} our
lifting conditions in the case of a complete Poisson group.
\section{Lifting of Lie algebroid homotopies}\label{llah}
Recall from \cite{crfs1}, that, for any Lie algebroid $A\rightarrow
M$ with anchor $\rho$,  an $A$-\tsf{path} is a  $\calC^1$ map
$\alpha:I\rightarrow A$, over a $\calC^2$ base path
$\gamma:I\rightarrow M$, such that  $\dd \gamma=(\rho\comp
\gamma)\,\alpha$, that is, a morphism of Lie algebroids 
$TI\rightarrow A$. For any Lie groupoid $\calG$, a
$\calG$-\tsf{path}  is a $\calC^2$ path within a source fibre,
starting from the unit section. Whenever $\calG$ is source connected
with Lie algebroid $A$, the right derivative 
$\delta_r:\{\calG\hbox{$-$}paths\}\to\{A\hbox{$-$}paths\}$, 
$\delta_r g(u)=\dd r_{g(u)}^{-1}\dot{g}(u)$, yields an homeomorphism
for the natural Banach topologies. Homotopy of $\calG$-paths within
the source fibres and  relative to the endpoints, called for short
$\calG$-\tsf{homotopy}, can be translated to an equivalence relation
for $A$-paths. An $A$-\tsf{homotopy} from an $A$-path $\alpha_-$ to
an $A$-path  $\alpha_+$ is a  $\calC^1$ morphism of Lie algebroids
$h:TI^{\hbox{\tiny{$\times 2$}}}\rightarrow A$, satisfying suitable
boundary conditions%
\footnote{ 
The boundary components of the compact square $I^{\times 2}$ shall be
denoted with  $\pa^\pm_{\mathrm{V}}=\{(\eps,1/2\pm1/2)\}_{\eps\in
I}$, respectively with  $\pa^\pm_{\mathrm{H}}=\{(1/2\pm1/2,u)\}_{u\in
I}$; $\iota_\pa:I\to I^{\times 2}$ is the inclusion of the boundary
component $\pa$.
}.  
Regarding $h$ as a 1-form taking values in the pullback $X^{\pmb{+}}
A$,  for the base map $X$, the anchor compatibility condition is $\dd
X=(\rho\comp X)\, h$,  while the bracket compatibility takes the form
of a Maurer-Cartan equation
$$
\qquad
\mathrm{D}\,h+\frac{1}{2}[h\,\overset{\wedge},\,h]=0
\qquad,\qquad 
\iota^*_{\pa^\pm_{\mathrm{H}}}h
=
\alpha_\pm\qquad\hbox{and}\qquad \iota^*_{\pa^\pm_{\mathrm{V}}}h
=
0
$$
where $\mathrm{D}$ is the covariant derivative of the pullback of an
arbitrary connection $\nabla$ for $A$ and $ [\theta_+\!
\,\overset{\wedge},\, \theta_-]$ is the contraction
$-\iota_{\theta_+\wedge\theta_-}X^{\pmb{+}}\tau^\nabla$ with the
pullback of the torsion tensor, $\theta_\pm\in\Ohm^1(X^{\pmb{+}} A)$.
It turns out that the quotients
$\{A\hbox{-}paths\}/A\hbox{-}homotopy$ and 
$\{\calG\hbox{-}paths\}/\calG\hbox{-}homotopy$ are both endowed with
natural Lie groupoid structures, isomorphic to the source 1-connected
cover of $\calG$, the groupoid multiplications essentially  being
given by concatenation in both cases.
\smallskip

A morphism $\phi:A\to B$ of Lie algebroids has the \cite{st}%
\smallskip\\
$l_0$) \tsf{0-\tla-homotopy lifting property} if, for any $A$-path
$\alpha_-$ and $B$-path $\beta_+$,  which is $B$-homotopic to
$\beta_-:=\phi\comp\alpha_-$,  there exists an $A$-path $\alpha_+$,
which is $A$-homotopic to $\alpha_-$ and satisfies 
$\phi\comp\alpha_+=\beta_+$;
\smallskip\\ 
$l_1$) \tsf{1-$\mathcal{LA}$-homotopy lifting property}  if,  for any
$A$-path $\alpha$, which is $A$-homotopic to the constant $A$-path
$\alpha_o\equiv 0_{\pr_A(\alpha(0))}$,  and $B$-homotopy $h_B$ from 
$\beta:=\phi\comp\alpha$ to the constant $B$-path  $\beta_o\equiv
0_{\pr_B(\beta(0))}$, there exists an $A$-homotopy $h_A$ from
$\alpha$ to $\alpha_o$,  such that $\phi\comp h_A=h_B$.
\smallskip\\ 
The \tla-homotopy lifting conditions above on a morphism of
integrable Lie algebroids  $A\to B$, translate to the infinitesimal
level homotopy lifting conditions in the  source fibres of the source
1-connected Lie groupoids $\calA$ and $\calB$ of $A$ and $B$ along
the integration $\calA\to\calB$.
\begin{proposition}\label{fib1} Let $A^{1,2}\to M^{1,2}$, $B\to N$ be
integrable Lie algebroids and $\phi_{1,2}:A^{1,2}\to B$ be
transversal morphisms of Lie  algebroids over $f_{1,2}:M^{1,2}\to N$,
so that the fibred product Lie algebroid $A^1\fib{\phi_1}{\phi_2}A^2$
is defined; denote with  $\fii_{1,2}:\calG^{1,2}\to \calH$ the
integrating morphisms for the  source 1-connected integrations. Then
the fibred product Lie groupoid
$\poidd{\calG^1\fib{\fii_1}{\fii_2}\calG^2}{M_1\fib{f_1}{f_2}M_2}$ 
is source 1-connected provided either $\phi_1$ or $\phi_2$ has the 0-
and 1- \tla-homotopy lifting property.
\end{proposition}
\noindent
This fact was proved in a (not so) special case in \cite{st}
(proposition 3.3.);  the same proof applies here, with the obvious
modifications. The transversality conditions are met, for example,
when either $\phi_1$ or $\phi_2$ is fibrewise submersive and base
submersive. 
\section{Dressing actions and symplectic double groupoids}\label{dasdg}
Let $(G,\Pi)$ be a Poisson group. The (\tsf{left}) \tsf{dressing
action} of $\frag^*$ on $G$ is the infinitesimal action
$\Upsilon:\frag^*\to\frax(G)$, $\Upsilon(\xi):=\Pi^\sharp\linv{\xi}$,
where $\linv{\xi}$ is the left invariant 1-form on $G$ associated
with $\xi\in\frag^*$. The (\tsf{left}) \tsf{dressing vector fields}
on $G$ are those in the image of $\frag^*$ under the left dressing
action map. Left dressing vector fields do  not have, in general,
complete flows; when they have, e.g. when $G$ is compact, $G$ is
called a \tsf{complete} Poisson group.
\bgn{remark} To any infinitesimal action $\sigma:\frag\to\frax(M)$
one can associate a Lie algebroid $\frag\acts M$ over $M$, whose
leaves are precisely the orbits of the infinitesimal action. The
anchor $\rho_\acts:\frag\acts M\to TM$ is given by
$\rho_\acts(x,q):=\sigma(x)_q$. Sections of $\frag\acts M$ are to be
identified with $\frag$-valued functions over $M$, thus
post-composition with the action map $\sigma$ induces a map
$\cif(M,\frag)\to\sf{End_{lin}}\cif(M,\frag)$,  $f\mapsto\sigma^f$, 
\be
\qquad
\pair{\sigma^{f_+}_q(f_-)}{\xi}
&:=\:&
(\sigma(f_+(q))_q(\pair{\xi}{f_-})\\
&\:=\:&
\rho_\acts(f_-)_q(\pair{\xi}{f_+})
\qquad,\qquad 
\xi\in\frag^*
\qquad,
\ee
$f_\pm\in\cif(M,\frag)$. The skewsymmetric bilinear operation
$$
\qquad
\bracts{f_+}{f_-}_q
:=
\brast{f_+(q)}{f_-(q)}
+
\sigma^{f_+}_q(f_-)
-
\sigma^{f_-}_q(f_+)
\quad,\quad 
q\in M
\quad,
$$	  
defines a Lie algebroid bracket on $\frag\acts M$ compatible with
$\rho_\acts$. Similarly, to any Lie group action $G\to\sf{Diff}(M)$,
one can associate an action Lie groupoid $\poidd{G\acts M}{M}$, whose
orbits on $M$ are the orbits of the action (we shall define action
groupoids for a Lie group(oid) actions in \S\ref{ipgs}).
\end{remark} 
In the case of the dressing action associated with a Poisson group
$G$, the left trivialization  $T^*G\to\frag^*\times G$ is an
isomorphism of Lie algebroids  to the action Lie algebroid.  When $G$
is complete, the infinitesimal dressing action integrates to a global
action, $\wt{\Upsilon}:G^\star\times G\to G$ of the dual 1-connected
Poisson group $G^\star$ on $G$ and  the action groupoid
$\poidd{G^\star\acts G}{G}$ is the source 1-connected integration of
$\frag\acts G$. Provided $G$ is 1-connected (in this case $G^\star$
is also complete)  the same argument applies to the integration of
the infinitesimal right dressing action  $\frag^*\to\frax(G)$,
defined using $\Pi^\star$, to a global action  $G^\star\times G\to
G^\star$.  One can show that the 1-connected integration $D$ of  the
Drinfel'd double $\mathfrak{d}$ is isomorphic to bitwisted product 
${G^\star{\ltimes\hs{-0.3}\rtimes} G}$ carrying two further
symplectic groupoid structures over $G$ and $G^\star$ \cite{ka,luth}
and making  $(G^\star{\ltimes\hs{-0.3}\rtimes} G, G ; G^\star
,\bullet)$ a double of the Poisson group on $G$. A construction of a
double in the  non-complete case, for a 1-connected $G$ was given in
\cite{lw} by Lu and Weinstein. We  remark that the total space of
Lu-Weinstein's double is  only locally diffeomorphic  to $D$, unless
$G$ is complete, and it is in general neither source
(1-)connected over $G$, nor over $G^\star$.  Moreover, the
1-connectivity of $G$ is essential in both constructions. 
The rest of this section is devoted to improve the above results to a
suitably simply connected integration, in the complete case,
dropping  the connectivity assumptions on $G$.
\bgn{theorem}\label{main} For any complete Poisson group $G$, the
source  1-connected symplectic groupoid $\calS$ of $G$ carries a
unique Lie groupoid structure over the 1-connected dual Poisson group
$G^\star$ making it a symplectic groupoid for the dual Poisson
structure and for which 
\bgn{equation}\label{mydouble}
\xy
*+{}="0",    <-0.7cm,0.7cm>
*+{\calS}="1", <0.7cm,0.7cm>
*+{G^\star}="2", <-0.7cm,-0.7cm>
*+{G}="3", <0.7cm,-0.7cm>
*+{\bullet}="4",
\ar  @ <0.07cm>   @{->} "1";"2"^{} 
\ar  @ <-0.07cm>  @{->} "1";"2"_{}  
\ar  @ <-0.07cm>  @{->} "1";"3"_{}
\ar  @ <0.07cm>   @{->} "1";"3"^{}
\ar  @ <0.07cm>   @{->} "2";"4"^{}
\ar  @ <-0.07cm>  @{->} "2";"4"_{}
\ar  @ <-0.07cm>  @{->} "3";"4"_{}
\ar  @ <0.07cm>   @{->} "3";"4"^{}
\endxy
\end{equation}
is a symplectic double groupoid.
\end{theorem}
\noindent
The symplectic double groupoid of last theorem is then a vertically
source 1-connected double of $G$. We derive this result from the
functorial approach to the integrability of $\lagpd$s developed in
\cite{st}; we shall recall below our approach in the special case of
the \tla-groupoid (fig. 3) of a Poisson group.	
\\
Let $\poidd{\calS}{G}$ be the source 1-connected symplectic groupoid
of $G$. The groupoid structural maps $\sh$, $\th$, $\epsh$ and
$\iotah$ of the cotangent prolongation groupoid
$\poidd{T^*G}{\frag^*}$ integrate uniquely to morphisms of Lie
groupoids $\nsth,\ntth:\calS\rightarrow G^\star$,
$\eth:G^\star\to\calS$ and $\ith:\calS\to\calS$. The nerves
$$
\qquad
(T^*G)^{\,\hbox{\tiny{$(n\!\!+\!\!1)$}}}
=
T^*G\fib{\sh}{\th\comp\pr_1}(T^*G)^{\,\hbox{\tiny{$(n)$}}}\qquad,\qquad
n\in\NN\qquad,
$$
of the top groupoid in figure 2, are canonically endowed with fibred
product Lie  algebroids over the side nerves 
$G^{\,\hbox{\tiny{$(n)$}}}=G^{\,\hbox{\tiny{$\times n$}}}$ and
integrate to the nerves of the top horizontal differentiable graph 
$(\calS,G^\star; \nsth,\ntth)$
$$
\qquad
\calS^{\,\hbox{\tiny{$(n\!\!+\!\!1)$}}}_{\hbox{\tiny{$H$}}}
=
\calS
\fib{\sth}{\tth\comp\pr_1}
\calS^{\,\hbox{\tiny{$(n)$}}}_{\hbox{\tiny{$H$}}}
\qquad,\qquad
n\in\NN\qquad.
$$
\emph{Provided the nerves
$\calS^{\,\hbox{\tiny{$(\bullet)$}}}_{\hbox{\tiny{$H$}}}$ are source
1-connected}, the top groupoid multiplication 
$\muh:(T^*G)^{\,\hbox{\tiny{$(2)$}}}\rightarrow T^*G$, integrates
uniquely to a morphism of Lie groupoids
$
\mth:\calS^{\,\hbox{\tiny{$(2)$}}}_{\hbox{\tiny{$H$}}}\to\calS
$.
One can check that it is indeed a groupoid multiplication compatible
with  $(\nsth,\ntth,\eth,\ith)$ by integrating the diagrams for the
groupoid structure on  $\poidd{T^*G}{\frag^*}$. In particular, it is
possible to integrate the associativity diagram, since the third
nerve $\calS^{\,\hbox{\tiny{$(3)$}}}_{\hbox{\tiny{$H$}}}$ is source
1-connected. Moreover, the graph $\gr{\muh}$ carries a Lie algebroid
structure over $\gr{\mu}$ making it a Lagrangian subalgebroid of
$T^*G\times T^*G\times\ol{T^*G}$; 
$\gr{\mth}\simeq\calS^{\,\hbox{\tiny{$(2)$}}}_{\,\hbox{\tiny{$H$}}}$
is the corresponding  source 1-connected integration, therefore
\cite{cat} a Lagrangian subgroupoid of 
$\calS\times\calS\times\ol{\calS}$. Thus (\ref{mydouble}) is indeed a
symplectic double groupoid and the induced Poisson structure on
$G^\star$ coincides  with $\Pi^\star$.
To complete the proof of theorem \ref{main}, we need to show that the
nerves $\calS^{\,\hbox{\tiny{$(\bullet)$}}}_{\hbox{\tiny{$H$}}}$ are
actually source 1-connected. This is consequence of proposition
\ref{fib1} and the following 
\bgn{lemma}\label{lift} For any complete Poisson group $G$, the
cotangent source map  $\sh:T^*G\to \frag^*$ has the 0- and
1-\tla-homotopy lifting properties.
\end{lemma}
\bgn{proof} Note that the diagram
$$
\xy
*+{T^*G}="0",    <2cm,0cm>
*+{\frag^*\acts G}="1", <1cm,-1cm>
*+{\frag^*}="2"
\ar  @ {->} "0";"1"^{\thicksim\quad}  
\ar  @ {->} "1";"2"^{\pr}
\ar  @ {->} "0";"2"_{\sh}
\endxy
$$
commutes in the category of Lie algebroids for the left
trivialization on the top edge: it suffices to prove the statement
for the projection   $\frag^*\acts G\to\frag^*$.%
\\ 
$\frag^*\acts G$-paths are pairs $(\xi,\gamma)$ of $\frag^*$-paths
and paths in $G$, such that
$\Upsilon(\xi(u))_{\gamma(u)}=\dd\gamma_u$, $u\in I$. A Lie algebroid
homotopy  $h_\acts: TI^{\hbox{\tiny{$\times 2$}}}\to\frag^*\acts G$
is a pair $(h, X)$ for which $h$ is a $\frag^*$-homotopy and $X:I\to
G$  satisfies  $\dd X=\Upsilon\comp h$. To see this, note that
$h_\acts$ takes values in the pullback  bundle  $X^+(\frag^*\acts
G)\simeq \frag^*\times I^{\hbox{\tiny{$\times 2$}}}$.  The bracket
compatibility for  $h_\acts$ can be written choosing the de Rham 
differential on $G$  with coefficients in $\frag^*$ as a linear
connection $\nabla$ for   $\frag^*\acts G$; this way the covariant
derivative for the pullback connection is simply the de Rham
differential on $I^{\hbox{\tiny{$\times 2$}}}$  with coefficients in
$\frag^*$. The torsion tensor of $\nabla$ is	
\be
\tau^\nabla(f_+,f_-)
&=&
{\dd f_-}_g(\rho^\acts(f_+))
-
{\dd f_+}_g(\rho^\acts(f_-))
-
\bracts{f_+}{f_-}_g\\
&=&
-\brast{f_+(g)}{f_-(g)}
\quad,
\ee
$f_\pm\in\cif(G,\frag^*)$, thus the pullback of $\tau^\nabla$ induces
the  canonical graded Lie bracket, also denoted by $\brast{\:}{\,}$,
on  $\Ohm^\bullet(I^{\hbox{\tiny{$\times 2$}}},\frag^*)$ and the
bracket compatibility condition for $h_\acts$ reduces to the
classical  Maurer-Cartan equation
$$
\dd h+\frac{1}{2}\brast{\,h\,}{h\,}=0
$$
for $h\in\Ohm^1(I^{\hbox{\tiny{$\times 2$}}},\frag^*)$. Suppose now
$(\xi_-,\gamma_-)$ is a fixed $\frag^*\acts G$-path, let $h$ be a
$\frag^*$-homotopy form $\xi_-$ to some other $\frag^*$-path $\xi_+$
and $H:I^{\hbox{\tiny{$\times 2$}}}\to G^\star$ the unique
$G^\star$-homotopy integrating $h$, i.e. such that 
$$
h=\delta_r^\eps H\cdot d\eps + \delta_r^uH \cdot du 
$$
for the partial right derivatives. We claim that $h_\acts:=(h,X)$,
where
\be
X(\eps,u)&=&H(\eps,u)\cdot H(0,u)\inverse *\gamma_-(u)\\
         &=&H(\eps,u)* H(0,u)\inverse *\gamma_-(u)
\qquad,\qquad u\:,\:\eps\in I\:\:,
\ee
is a $\frag^*\acts G$-homotopy; here the symbol $*$ denotes the
integrated  dressing action. We postpone to the appendix this
straightforward but lengthy check. The lifting conditions follow: by
construction
\be
\iota^*_{\pa^\pm_V}(h,X)&=&(\iota^*_{\pa^\pm_V}h,X\comp\iota_{\pa^\pm_V})=0\\
\iota^*_{\pa^\pm_H}(h,X)&=&(\xi_\pm, \gamma_\pm) 
\ee
for some path $\gamma_+$ in $G$, 
($l_0$) If $h$ is a fixed $\frag^*$-homotopy to some fixed 
$\frag^*$-path $\xi_+$, $(\xi_+,\gamma_+)$ is the desired lift. 
($l_1$) In particular if $\xi_-$ is $\frag^*\acts G$-homotopic to the
constant  $\frag^*\acts G$-path $\xi_o\equiv 0$ and $h$ a
$\frag^*$-homotopy from $\xi_-$ to the constant $\frag^*$-path
$\xi_+\equiv 0$, we have
$$
\qquad
\dd\gamma_+|_u=\Upsilon(\xi_+(u))_{\gamma_+(u)}=0_{\gamma_+(u)}
\qquad, 
$$
hence $ \gamma_+(u)\equiv\gamma_+(0)=\gamma_-(0)$, since the base
paths of homotopic Lie algebroid paths are homotopic relatively to
the endpoints; $(h,X)$ is then  the desired homotopy.
\end{proof}
\section{Integrability of Poisson $G$-spaces}\label{ipgs}
A \tsf{Lie groupoid action} of a  Lie groupoid $\gpdm$ on a smooth
map $j:P\to M$ is given by a smooth map
$\sigma:\calG\fib{\sor}{j}P\to P$, $(g,p)\mapsto g*p$, satisfying the
conditions 
\bgn{equation}\label{action}
j(g*p)=\tar(g)
\qquad
\eps(j(n))*n=n
\qquad
g*(h*p)=(g\cdot h)*p
\qquad,
\end{equation}
to be understood whenever they make sense, that generalize the notion
of a  Lie group action. The map $j$ above is called \tsf{moment map}
of the action. To any Lie groupoid action one can associate an action
Lie groupoid  $\poidd{\calG\acts P}{P}$ with total space
$\calG\fib{\sor}{j}P$, whose the source map is the first projection,
the target map coincides with the action map, the partial
multiplication is given by
$$
\qquad(g_+,p_+)\cdot_\acts(g_-,p_-)=(g_+\cdot g_-, p_-)\qquad,
$$
for all elements of $\calG$ and $P$ for which the expression above
makes sense; unit section and inversion are defined accordingly.%
\\
When $\calG$ is a Poisson groupoid, $ P$ a Poisson manifold and $j$ a
Poisson map for the induced Poisson structure on $M$,  $P$ is called
a \tsf{Poisson $\calG$-space} if  
\bgn{equation}\label{compa}
\gr{\sigma}\subset\calG\times P\times\ol{P}
\end{equation}
is a coisotropic submanifold; in that case one calls the action  a
\tsf{Poisson action}. The same definition applies when $\calG$ is a
Poisson group or a  symplectic groupoid acting on a symplectic
manifold. The compatibility condition is equivalent to require the 
action map $G\times P\to P$ to be Poisson, in first case, and to
$\gr{\sigma}$ being Lagrangian, in the second, this makes the action
a \emph{symplectic groupoid action} in the sense of Mikami and
Weinstein \cite{MiWs88}.  
Let $P$ be a Poisson $G$-space; dualizing the infinitesimal action
$\sigma_o:\frag\to\frax(P)$, yields a morphism of Lie bialgebroids
$T^*P\to \frag^*$,
\bgn{equation}\label{moment}
\qquad
\pair{\jh(\alpha_p)}{x}:=\pair{\alpha_p}{\sigma_o(x)}\qquad,
\qquad x\in\frag\qquad,
\end{equation}
i.e. $\jh$ is at the same time a morphism of Lie algebroids and a
Poisson map for the dual Poisson structures \cite{xupg}. Remarkably,
the $G$-action cotangent lifts to a Poisson action of the symplectic
groupoid $\poidd{T^*G}{\frag^*}$ on the moment map $\jh$; moreover
the lifted action is compatible with the Lie algebroids on $T^*G$ and
$T^*P$ \cite{mu}. In fact, the lifted action map
$\hat{\sigma}:T^*G\fib{\sh}{\jh}T^*P\to T^*P$, 
$
\theta_g*\alpha_p
:=
\hat{\sigma}(\theta_g,\alpha_p):=\alpha_p\comp\dd\sigma_{g^{-1}}
$, 
is a morphism of Lie algebroids and it follows that
\bgn{equation}\label{lapga}
\xy
*+{}="0",    <-1.2cm,0.7cm>
*+{T^*G\acts T^*P}="1", <0.9cm,0.7cm>
*+{T^*P}="2", <-1.2cm,-0.7cm>
*+{G\acts  P}="3", <0.9cm,-0.7cm>
*+{P}="4",
\ar  @ <-0.07cm>   @{->} "1";"2"^{} 
\ar  @ <0.07cm>    @{->} "1";"2"^{}  
\ar  		   @{->} "1";"3"^{}
\ar                @{->} "2";"4"_{}
\ar  @ <0.07cm>    @{->} "3";"4"^{}
\ar  @ <-0.07cm>   @{->} "3";"4"_{}
\endxy
\end{equation}
is an $\lagpd$, for the action groupoid of the cotangent lifted
action on the top horizontal side and the fibred product Lie
algebroid on the top vertical side.%
\\
Mackenzie's construction of (\ref{lapga}) should be regarded as the
lift of a Poisson group action to a morphic action  in the category
of Lie algebroids. By analogy, it is natural to consider morphic
actions in the category of Lie groupoids, that is,  pairs of groupoid
actions  $(\widetilde{\sigma}, \sigma)$ on pairs  of moment maps
$(\calJ,j)$, such that both $(\calJ,j)$ and 
$(\widetilde{\sigma},\sigma)$ form morphisms of Lie groupoids (fig.
5, 6). Thanks to the regularity condition on the double source map of
$\calD$, the fibred product
$\poidd{\calD\fib{\sth}{\calJ}\calB}{\calH\fib{\ssh}{j}N}$  is always
a Lie groupoid; moreover the diagram in figure 7 is a double Lie
groupoid, we shall  refer to as the double Lie groupoid of the
action, for the action groupoids on the horizontal edges.
$$
\bgn{array}{ccccc}
{\xy
*+{}="0",    <-0.9cm,0.7cm>
*+{\calD\fib{\sth}{\calJ}\calB}="1", <0.9cm,0.7cm>
*+{\calB}="2", <-0.9cm,-0.7cm>
*+{\calH\fib{\ssh}{j}N}="3", <0.9cm,-0.7cm>
*+{N}="4",
\ar  @ {->} "1";"2"^{\qquad\Tilde{\sigma}}  
\ar  @ <-0.07cm>  @{->} "1";"3"_{}
\ar  @ <0.07cm>   @{->} "1";"3"^{}
\ar  @ <0.07cm>   @{->} "2";"4"^{}
\ar  @ <-0.07cm>  @{->} "2";"4"_{}
\ar  @ {->} "3";"4"^{\qquad\sigma}
\endxy}
&\qquad\qquad&
{\xy
*+{}="0",    <-0.7cm,0.7cm>
*+{\calB}="1", <0.7cm,0.7cm>
*+{\calV}="2", <-0.7cm,-0.7cm>
*+{N}="3", <0.7cm,-0.7cm>
*+{M}="4",
\ar  @ {->} "1";"2"^{\calJ}  
\ar  @ <-0.07cm>  @{->} "1";"3"_{}
\ar  @ <0.07cm>   @{->} "1";"3"^{}
\ar  @ <0.07cm>   @{->} "2";"4"^{}
\ar  @ <-0.07cm>  @{->} "2";"4"_{}
\ar  @ {->} "3";"4"^{j}
\endxy}
&
\qquad
\qquad
&
{\xy
*+{}="0",    <-0.9cm,0.7cm>
*+{\calD\acts\calB}="1", <0.9cm,0.7cm>
*+{\calB}="2", <-0.9cm,-0.7cm>
*+{\calH\acts N}="3", <0.9cm,-0.7cm>
*+{N}="4",
\ar  @ <0.07cm>    @ {->} "1";"2"
\ar  @ <-0.07cm>   @ {->} "1";"2"
\ar  @ <-0.07cm>   @ {->} "1";"3"_{}
\ar  @ <0.07cm>    @ {->} "1";"3"^{}
\ar  @ <0.07cm>    @ {->} "2";"4"^{}
\ar  @ <-0.07cm>   @ {->} "2";"4"_{}
\ar  @ <-0.07cm>   @ {->} "3";"4"
\ar  @ <0.07cm>    @ {->} "3";"4"
\endxy}\\
Figure\: 5
&\qquad\qquad&
Figure\:6
&
\qquad
\qquad
&
Figure\:7
\end{array}
$$
It is clear that a morphic action in the category of Lie groupoids
differentiates to a morphic action in the category of Lie algebroids.
A natural question is then whether the cotangent lift of a Poisson
$G$-space integrates to a morphic action of a symplectic double
groupoid integrating $G$ on the integration $\calJ:\Lambda\to
G^\star$ of (\ref{moment}). The answer is positive in the complete
case.
\bgn{theorem}\label{mainapp} Let $G$ be a complete Poisson group and
$P$ an integrable Poisson manifold with source 1-connected symplectic
groupoid $\Lambda$.  If $P$ is a Poisson $G$-space, then
\smallskip\\
$i)$ $\calJ:\Lambda\to G^\star$ is a Poisson $\calS$-space for the
top horizontal groupoid $\poidd{\calS}{G^\star}$ of the vertically
source 1-connected double of $G$;
\smallskip\\
$ii)$ The action map 
$
\wt{\sigma}:\calS\fib{\sth}{\calJ}\Lambda\to\Lambda
$ 
is a morphism of Lie groupoids over the action map $\sigma :G\times
P\to P$;
\smallskip\\
$iii)$ The action double groupoid of the integrated action
\bgn{equation}\label{lgsga}
\xy
*+{}="0",    <-0.7cm,0.7cm>
*+{\calS\acts \Lambda}="1", <0.7cm,0.7cm>
*+{\Lambda}="2", <-0.7cm,-0.7cm>
*+{G\acts P}="3", <0.7cm,-0.7cm>
*+{P}="4",
\ar  @ <0.07cm>   @{->} "1";"2"^{} 
\ar  @ <-0.07cm>  @{->} "1";"2"_{}  
\ar  @ <-0.07cm>  @{->} "1";"3"_{}
\ar  @ <0.07cm>   @{->} "1";"3"^{}
\ar  @ <0.07cm>   @{->} "2";"4"^{}
\ar  @ <-0.07cm>  @{->} "2";"4"_{}
\ar  @ <-0.07cm>  @{->} "3";"4"_{}
\ar  @ <0.07cm>   @{->} "3";"4"^{}
\endxy
\end{equation}
is the vertically source 1-connected double Lie groupoid integrating
the $\lagpd$ {\em (\ref{lapga})} associated with the Poisson action.
\end{theorem} 
\bgn{proof} It was proved in \cite{xupg} that $\calJ$ is Poisson. We
claim that the  fibred product Lie groupoids
$\calS\fib{\sth}{\calJ}\Lambda$ and 
$$
\calS\fib{\sth}{\calJ}(\calS\fib{\sth}{\calJ}\Lambda) \simeq
(\calS\fib{\sth}{\tth}\calS)\fib{\sth\comp \pr_2}{\calJ}\Lambda 
$$
are source 1-connected; then $\sigmah\comp(\id_{T^*G}\times
\sigmah)$, and $\sigmah\comp(\muh\times\id_{T^*P})$ integrate to the
morphisms  $\wt{\sigma}\comp(\id_{\calS}\times\wt{\sigma})$ and
$\wt{\sigma}\comp(\mth\times\id_{\Lambda})$, which coincide up to the
identification of the domains, by uniqueness. Compatibility of
$\wt{\sigma}$ with $\calJ$, follows in the same fashion. 
($ii$) holds by construction and ($iii$) is now obvious. 
($i$) $\gr{\widetilde{\sigma}}\simeq\calS\fib{\sth}{\calJ}\Lambda$ 
is the source 1-connected integration of 
$$
\qquad
\gr{\hat{\sigma}}
=
(\sf{id}_{T^*G}\times\sf{id}_{T^*P}\times-\sf{id}_{T^*P})\,N^*\gr{\sigma}
\qquad,
$$
which is a Lagrangian Lie subalgebroid of $T^*G\times T^*P\times
\ol{T^*P}$, by coisotropicity of $\gr{\sigma}$ in $\calG\times
P\times \ol{P}$; it follows \cite{cat} that
$\gr{\widetilde{\sigma}}\subset \calS\times\Lambda\times\ol{\Lambda}$
is a Lagrangian subgroupoid.
Our claim follows from proposition \ref{fib1} if the projection
$T^*G\acts T^*P\to T^*P$ has the lifting properties ($l_0$, $l_1$).
Note that  $T^*G\acts T^*P$-paths, -homotopies, are pairs
$(X^{\sg},X^{\hbox{\tiny{$P$}}})$ of $T^*G$-paths, -homotopies, and
$T^*P$-paths, -homotopies, such that $\sh\comp X^{\sg}=\jh\comp
X^{\hbox{\tiny{$P$}}}$.
For any assigned ($T^*G\acts T^*P$)-path
$(\xi^{\sg}_-,\xi^{\hbox{\tiny{$P$}}}_-)$ and  $T^*P$-homotopy $h^P$
from $\xi^{\hbox{\tiny{$P$}}}_-$ to some $T^*P$-path 
$\xi^{\hbox{\tiny{$P$}}}_+$,  $\jh\comp h^P$ is a $\frag^*$-homotopy
from $\sh\comp\xi^{\sg}_-$ to some $\frag^*$-path $\xi_+$; hence by
applying lemma \ref{lift} $\jh\comp h^P$ can be lifted to a
$T^*G$-homotopy $h^{\sg}$ such that $\sh\comp h^{\sg}=\jh\comp
h^{\hbox{\tiny{$P$}}}$ and $\pmh h^{\sg}=\xi_-^{\sg}$. Thus the pair
$H^\acts:= (h^{\sg},h^{\hbox{\tiny{$P$}}})$ is a  $T^*G\acts
T^*P$-homotopy lifting $h^{\hbox{\tiny{$P$}}}$, such that $\pmh
H^\acts=(\xi^{\sg}_-,\xi^{\hbox{\tiny{$P$}}}_-)$; as in the proof of
lemma \ref{lift},  the boundary conditions necessary to fulfill
($l_0$, $l_1$) hold according to the initial data.
\end{proof}
\bgn{remark} For any symplectic groupoid action such as above,  the
multiplicativity condition 
\bgn{equation}\label{multisigma}
\wt{\sigma}^*\ohm=\pr_\calS^*\Ohm+\pr_\Lambda^*\ohm
\end{equation}
on  $\calS\acts\Lambda$ is equivalent to  $\gr{\wt{\sigma}}\subset
\calS\times\Lambda\times\ol{\Lambda}$ being Lagrangian, for the
symplectic forms $\Ohm$ on $\calS$ and $\ohm$ on $\Lambda$.
Explicitly, (\ref{multisigma}) reads 
$$
\qquad
\ohm_{s*\lambda}(\delta s_+*\delta \lambda_+,\delta s_+*\delta \lambda_+) =
\Ohm_{s}(\delta s_+,\delta s_+) + \ohm_{\lambda}(\delta \lambda_+,\delta \lambda_+)
\qquad, 
$$
where we have used the symbol $*$ for the tangent lift of
$\wt{\sigma}$, for all composable  $\delta s_\pm\in T\calS$ and
$\delta \lambda_\pm\in T\Lambda$.
\end{remark}
\section{Integrability of quotient Poisson structures}\label{iqps}
Consider any Poisson action $\sigma:G\times P\to P$ of a (not
necessarily complete) Poisson group $G$. By coisotropicity of
$\gr{\sigma}$  the $G$-invariant functions on P form a Poisson
subalgebra  $\cif(P)^{\sg}\subset\cif(P)$. When the quotient $P/G$ is
smooth, i.e. if the action is free and proper,  its algebra of
functions coincides with $\cif(P)^{\sg}$ and  can be endowed with a
Poisson bracket simply by setting 
$$
\qquad
\poib{f_+}{f_-}_{P/G}([p]_{P/G}) = \poib{f_+}{f_-}_P(p)\qquad,\qquad
f_\pm\in\cif(P)^{\sg}\qquad,
$$
for any choice of a representative $p\in[p]_{P/G}\in P/G$. Therefore
we have
\bgn{proposition}\cite{wsc}\label{known} Let a Poisson group $G$ act
on a Poisson manifold $P$ freely and properly. If the action is
Poisson, there exists a unique Poisson structure on $P/G$ making the
quotient projection a Poisson submersion. 
\end{proposition}
Next, we shall obtain a symplectic groupoid for the quotient Poisson
manifold, when $P$ is integrable and $G$ complete, using a reduction
procedure on the action double groupoid of theorem \ref{mainapp}; a
few remarks are in order.
\\
On the one hand, by the very definition (\ref{moment}) of the moment
map $\jh$, $\ker_p\jh=N^*_p\calO$, for the $G$-orbit  $\calO$ through
$p\in P$; thus, whenever the action is free,
$$
\rank\jh=\dim P-(\dim P -\dim \calO)=\dim\frag^*
$$
and $\jh$ is a bundle map of maximal rank. It follows, if $P$ is
integrable, that the integrating morphism $\calJ:\Lambda\to G^\star$
is sourcewise submersive and submersive, therefore the kernel
groupoid $\ker \calJ=\calJ\inverse(e_\star)$ is a wide \emph{Lie}
subgroupoid.
$$
\xy
*++{}="0",    <-1.3cm,0.7cm>
*++{\calS\acts\Lambda}="1", <0.7cm,0.7cm>
*++{\Lambda}="2", <-1.3cm,-0.7cm>
*++{G\acts  P}="3", <0.7cm,-0.7cm>
*++{P}="4",     <2.2cm,-0.7cm>
*++{\calS}="1'", <4.2cm,-0.7cm>
*++{G^\star}="2'", <2.2cm,-2.1cm>
*++{G}="3'", <4.2cm,-2.1cm>
*++{\bullet}="4'",  <-1cm,-2.1cm>
*++{Figure\: 8}="2''"
\ar  @ <-0.07cm>   @{->} "1";"2"^{} 
\ar  @ <0.07cm>    @{->} "1";"2"^{}  
\ar   @ <0.07cm>   @{->} "1";"3"_{}
\ar   @ <-0.07cm>  @{->} "1";"3"_{}
\ar   @ <0.07cm>   @{->} "2";"4"_{}
\ar   @ <-0.07cm>  @{->} "2";"4"_{}
\ar  @ <0.07cm>    @{->} "3";"4"^{}
\ar  @ <-0.07cm>   @{->} "3";"4"_{}
\ar  @ <-0.07cm>   @{->} "1'";"2'"^{} 
\ar  @ <0.07cm>    @{->} "1'";"2'"^{}  
\ar   @ <0.07cm>   @{->} "1'";"3'"_{}
\ar   @ <-0.07cm>  @{->} "1'";"3'"_{}
\ar   @ <0.07cm>   @{->} "2'";"4'"_{}
\ar   @ <-0.07cm>  @{->} "2'";"4'"_{}
\ar  @ <0.07cm>    @{->} "3'";"4'"^{}
\ar  @ <-0.07cm>   @{->} "3'";"4'"_{}
\ar  		   @{->} "1";"1'"^{\hs{-1.7}\pr_\calS}
\ar  		   @{->} "2";"2'"^{\calJ}
\ar  		   @{->} "3";"3'"^{}
\ar  		   @{->} "4";"4'"^{j}
\endxy
\qquad\qquad
\xy
*+{}="0",    <-1cm,0.7cm>
*+{G\acts\hbox{\small{$\jj$}}}="1", <1.4cm,0.7cm>
*+{\hbox{\small{$\jj$}}}="2", <-1cm,-0.7cm>
*+{G\acts P}="3", <1.4cm,-0.7cm>
*+{P}="4", <0.2cm,-2.1cm>
*+{Figure\: 9}="sarcazzo"
\ar  @ <0.07cm>   @{->} "1";"2"^{} 
\ar  @ <-0.07cm>  @{->} "1";"2"_{}  
\ar  @ <-0.07cm>  @{->} "1";"3"_{}
\ar  @ <0.07cm>   @{->} "1";"3"^{}
\ar  @ <0.07cm>   @{->} "2";"4"^{}
\ar  @ <-0.07cm>  @{->} "2";"4"_{}
\ar  @ <-0.07cm>  @{->} "3";"4"_{}
\ar  @ <0.07cm>   @{->} "3";"4"^{}
\endxy
$$
On the other hand, the moment map for the integrated action can
always be completed to a morphism of double Lie groupoids (fig. 8),
whose (vertical) kernel (fig. 9) is  naturally an action double
groupoid for the restriction of the top action.
\bgn{theorem}\label{redux} Let $G$ be a complete Poisson group and
$P$ a Poisson $G$-space. If $G$ acts freely and properly and $P$ is
integrable, the quotient Poisson manifold $P/G$ is also integrable;
in particular
\smallskip\\
$i)$ The quotient $\jmod$ is smooth  and the unique symplectic form 
$\ul{\ohm}$  on $\jmod$  such that
\bgn{equation}\label{qsymp}
\qquad\pr^*\ul{\ohm}=\iota^*\ohm\qquad
\end{equation}
makes it a symplectic groupoid over $P/G$;
\smallskip\\ 
$ii)$ 
$\poidd{\jmod}{P/G}$ integrates the quotient Poisson structure.
\end{theorem}
\bgn{proof} Note that $\calK:=\jj$ is coisotropic since $\calJ$ is
Poisson and $e_\star\subset G^\star$ is coisotropic. ($i$) There is a
natural representation $\vartheta:G\rightarrow\sf{Aut}(\calK)$, 
$
\vartheta_g(\kappa)
=
\widetilde{\sigma}(\eps_{\hbox{\tiny{$V$}}}(g),\kappa)=g*k
$ 
in the group of Lie groupoid automorphisms, providing the descent
data for the Lie groupoid structure. Since $P\to P/G$ is a principal
$G$-bundle, it follows from (\cite{mrmc}, lemma 2.1) that $\calK/G$
is a Lie groupoid over $P/G$.
On $G\acts\calK$ the
multiplicativity condition reads
$$
\qquad
\ohm_{g*\kappa}(\delta g_+*\delta \kappa_+,\delta g_-*\delta \kappa_-)
=
\ohm_{\kappa}(\delta \kappa_+,\delta \kappa_-)
\qquad,
$$
being $G\subset\calS$ Lagrangian, for all composable  $\delta
g_\pm\in TG$ and $\delta \kappa_\pm\in T\calK$; it follows that 
$\vartheta_g^*(\iota^*\ohm)\equiv\iota^*\ohm$, hence using
(\ref{qsymp}) to define a closed 2-form on $\calK/G$ makes sense. We
claim that the characteristic distribution $\Delta$ of ${\calK}$
spans $T_\kappa\calO$ for the $\calG$-orbit $\calO$ through
$\kappa\in\calK$; it follows that $\ul{\ohm}$ is nondegenerate. To
see this, note that
\be
\rank T^\ohm\calK
&=&
\dim \Lambda -\rank \calJ =\dim G^\star\\
&=&\dim \calO
\ee
and that all $\delta o\in T_k\calO$ are of the form  $\delta o=\delta
g*0_k$ for  the tangent lifted action, thus, for all   $\delta k\in
T_k\calK$,   
$$
\qquad
\ohm_{k}(\delta k,\delta g*0_k)
=
\ohm_{k}(\delta g\inverse*\delta k, 0_k)=0\qquad,
$$ 
i.e. for all $ k\in\calK$
$$
\qquad
\Delta_k
=
T^\ohm_k\calK=T_k\calO\qquad.
$$
Multiplicativity of $\ul{\ohm}$ follows easily from multiplicativity
of $\ohm$ and uniqueness is clear thanks to condition (\ref{qsymp}). 
($ii$) Since the characteristic leaves of $\calK$ are the connected
components of the $G$-orbits, 
$\cif(\calK)^G\subset\cif(\calK)^\Delta=\ker\Delta$, thus all
extension $F_\pm\in\cif(\Lambda)\,$ of $\,f_\pm\in\cif(\calK)^G$ are
in the  normalizer of the vanishing ideal of $\,\calK\,$ and the
Hamiltonian vector fields $\,X^{F_\pm}$ restrict to $\,\calK\,$. 
\noindent
Let $X^{F_\pm}=\dd\iota Y^{F_\pm}$ for some (in general non smooth)
vector fields $Y^{F_\pm}$ on $\calK$. We have
\be
(\dd\pr^t_{\pr(\kappa)}\comp\ul{\ohm}^\sharp_{\pr(\kappa)}\comp\dd\pr_{k})Y^{F_\pm}
&=&
(\dd\iota^t_{\iota(\kappa)}\comp\ohm^\sharp_{\iota(\kappa)}) X^{F_\pm}_{\iota(\kappa)}\\
&=&\dd_{\calK}f_\pm\\
&=&
\dd\pr^t_{\pr(\kappa)}\dd_{\calK/G}f_\pm
\quad,
\ee
i.e. the Hamiltonian vector fields $\ul{X}^{f_\pm}$ of
$f_\pm\in\cif(\calK/G)$ are given by
$\ul{X}^{f_\pm}_{\pr(\kappa)}=\dd\pr Y^{F_\pm}_\kappa$ and the
Poisson bracket $\poib{}{}_{\calK/G}$ associated with $\ul{\ohm}$ can
be computed using extensions:
\be
\pr^*\poib{f_+}{f_-}_{\calK/G}
&:=\:&
\ul{\ohm}(\ul{X}^{f_-},\ul{X}^{f_+})\comp\pr\\
&\:=\:&
\pr^*\ul{\ohm}(Y^{F_-},Y^{F_+})\\
&\:=\:&
\ohm(X^{F_-},X^{F_+})\comp\iota\\
&\:=:&\iota^*\poib{F_+}{F_-}_{\Lambda}\qquad\qquad.
\ee
Let now $\poib{}{}'$ be the Poisson bracket induced by the symplectic
groupoid of ($i$)  on $P/G$ and $u_\pm\in\cif(P)^G$; since 
$\sor_\Lambda^*u_\pm\in\cif(\Lambda)$ extend
$\sor_{\calK}^*u_\pm\in\cif(\calK)^G$, 
\be
\poib{u_+}{u_-}'([p]_{P/G})
&:=\:&
\poib{\sor_{\calK/G}^*u_+}{\sor_{\calK/G}^*u_-}_{\calK/G}(\eps_{\calK/G}([p]_{P/G}))\\
&\:=\:&
\poib{\sor_\Lambda^*u_+}{\sor_\Lambda^*u_-}_{\Lambda}(\eps_\Lambda(p))\\
&\:=\:&
\poib{u_+}{u_-}_{P}(p)\\
&\:=:&
\poib{u_+}{u_-}_{P/G}([p]_{P/G})\hs{3},
\ee
for all $p\in P$, ($ii$) follows uniqueness (proposition \ref{known}).
\end{proof}
Theorem \ref{redux} generalizes a result by Xu (\cite{xu}, theorem
4.2), regarding Poisson actions with a complete momentum map. A
\tsf{momentum map} \cite{lu} for a Poisson $G$-space $(P,\pi)$ is a
Poisson map  $j:P\to G^\star$ such that $\sigma(x)=\pish
j^*\linv{x}$, for the infinitesimal action
$\sigma(\cdot):\frag\to\frax(P)$ and the left invariant 1-form
$\linv{x}$ on $G^\star$ associated with  $x\in\frag\simeq\frag^{**}$;
such a Poisson map is called complete, when the Hamiltonian vector
field of $j^*f$ is complete for all compactly supported 
$f\in\cif(G^\star)$. If $G$ is 1-connected, so that the right
dressing action of $G$ on $G^\star$ is globally defined, $j$ is
always equivariant. The construction of \cite{xu}, when $G$ is
complete and 1-connected, the action admits a complete momentum map 
$j$, produces a symplectic groupoid
$\poidd{J\inverse(e_\star)/G}{P/G}$, assuming the quotient 
$J\inverse(e_\star)/G$ is smooth, where the momentum map 
$J:\Lambda\to G^\star$, is given by $J(\lambda)=j(\tar(\lambda))\cdot
j(\sor(\lambda))\inverse$.  Note that $J$ is by construction a
morphism of Lie groupoids and one can check \cite{xupg} that it
differentiates to $\jh$, therefore $J$ coincides with our $\calJ$; it
is easy to see that the $G$ action on $J\inverse(e_\star)$ of
\cite{xu} is the same as that induced by the
$\poidd{\calS}{G^\star}$-action on $\Lambda$   ($\calS\simeq
G^\star{\ltimes\hs{-0.3}\rtimes} G$, under the assumptions).
Specializing last theorem to the case considered in \cite{xu}, shows
that Xu's quotient is always smooth. 
A Lie group $G$ is trivially ($\Pi=0$) a complete Poisson group, with
the abelian group $\ol{\frag^*}$ as a dual Poisson group; in this
case a Poisson group action is an action by Poisson diffeomorphisms
and our approach of sections  \S\ref{dasdg}-\S\ref{ipgs} reproduces
the ``symplectization functor'' treatment of Fernandes \cite{fsf} and
Fernandes-Ortega-Ratiu \cite{for} from the viewpoint of double
structures. 
The construction given in \cite{for} (proposition 4.6) of a
symplectic groupoid for  the quotient Poisson manifold is precisely
the construction of theorem  \ref{redux} in the special case $\Pi=0$.
In the same paper,  a condition (\cite{for}, theorem 4.11) for
integration to commute with reduction  is also given. This issue
depends only on the connectivity of the source fibres of 
$\poidd{\jmod}{P/G}$ and the conditions of \cite{for} apply to
Poisson actions of Poisson groups. From the proof of theorem
\ref{redux} it is clear that the Lie groupoid on  $\jmod$ is really a
``push-forward'' of $\jj$, obtained by identifying the source fibres
along the $G$-orbits of $P$; then  $\jmod$ is source (1-)connected
iff so is $\jj$.%
\\
For all $p\in P$, the quotient $\KK_p(\jh)$ of the space of
$T^*P$-loops taking values in $\ker\,\jh$ which are $T^*P$-homotopic
to the null $T^*P$-path $0_p$ modulo $T^*P$-homotopies taking values
in  $\ker\,\jh$, is  naturally a group.
\bgn{proposition} Under the same hypotheses of theorem \ref{redux},
the source connected component of $\poidd{\jmod}{P/G}$ is the source
1-connected symplectic groupoid of $P/G$ iff the groups $\KK_p(\jh)$
are trivial, for all $p\in P$.
\end{proposition}
The proof of (\cite{for}, theorem 4.11) applies for $\Pi\neq 0$ with
no modification.
\appendix
\section{End of proof of proposition \ref{lift}}
\noindent
It remains to check the anchor compatibility condition. Recall that
\be
X(\eps,u)&=&H(\eps,u)\cdot H(0,u)\inverse *\gamma_-(u)\\
         &=&H(\eps,u)* H(0,u)\inverse *\gamma_-(u)
\qquad,\qquad u\:,\:\eps\in I\:\:,
\ee
and set $a=\delta_r^uH$ and  $b=\delta_r^\eps H$. The derivative of
$X$ in the $\eps$-direction is thus
\bea
\pa_\eps X(\eps,u)
&=&
\Upsilon(\pa_\eps H(\eps,u))_{H(0,u)\inverse *\gamma_-(u)}
\:=\:
\Upsilon(\dd r_{H(\eps,u)}b(\eps,u))_{H(0,u)\inverse *\gamma_-(u)}\\
&=&
\Upsilon(b(\eps,u))_{X(\eps,u)}\hs{3},
\eea
since for all $g\in G$, $h^\star\in G^\star$,
$\xi=\dot{g}^\star(o)\in\frag^*$ and path $g^\star$ in $G^\star$
$$
\Upsilon(\dd r_{h^\star}\xi)
=
\left.
\frac{d}{d\alpha}\right|_{\alpha=o}
\wt{\Upsilon}(g^\star(\alpha)\cdot h^\star, g)
=
\left.
\frac{d}{d\alpha}\right|_{\alpha=o}
\wt{\Upsilon}(g^\star(\alpha),h^\star*g)
=
\Upsilon(\xi)_{h^\star*g}
\quad.
$$
The computation of the derivative in the $u$-direction is more
involved. We have
$$
\qquad
\pa_u X(\eps,u)
=
\dd\wt{\Upsilon}_{(H(\eps,u)\cdot H(0,u)\inverse, \gamma_-(u))}
		(\delta	 H, \dot{\gamma}_-(u))
\qquad,
$$
where 
$\quad
\dot{\gamma}_-(u)
=
\Upsilon(\xi_-(u))_{\gamma_-(u)}
=
\Upsilon((\iota_{\pa\hup^{-}}^* h)(u))
=
\Upsilon(a(0,u))_{X(0,u)}
\quad$ and
\bea
\delta H 
&=&
\dd\mu^\star_{(H(\eps,u),H(0,u)\inverse)}
(\pa_u H(\eps,u),\dd\iota_\star\pa_u H(0,u))\\
&=&
\dd\mu^\star_{(H(\eps,u),H(0,u)\inverse)}
(\dd r_{H(\eps,u)}a(\eps,u),\dd\iota_\star\dd r_{H(0,u)}a(0,u))\\
&=&
\dd\mu^\star_{(H(\eps,u),H(0,u)\inverse)}
(0_{H(\eps,u)},\dd\iota_\star\dd r_{H(0,u)}a(0,u))\\
&+&
\dd\mu^\star_{(H(\eps,u),H(0,u)\inverse)}
(\dd r_{H(\eps,u)}a(\eps,u),0_{H(0,u)\inverse})\\
&=& \delta H_+ +\delta H_-\hs{5},
\eea
since the tangent grop multiplication 
$\dd\mu^\star:TG^\star\times TG^\star\to TG^\star$ 
is fibrewise linear, with 
\bea
\delta H_+
&=&
\dd l_{H(\eps,u)}\dd\iota_\star\dd r_{H(0,u)}a(0,u)
\:\:=\:\:
\dd l_{H(\eps,u)}\dd l_{H(0,u)\inverse}\dd\iota_\star a(0,u)
\\
&=&
\dd l_{H(\eps,u)\cdot H(0,u)\inverse}\dd\iota_\star a(0,u)
\\
\delta H_-&=&\dd r_{H(0,u)}\inverse\dd r_{H(\eps,u)}a(\eps,u)=
\dd r_{H(\eps,u)\cdot H(0,u)\inverse}a(\eps,u)\hs{2}.
\eea
The tangent action map $\dd\wt{\Upsilon}:TG^\star\times TG\to TG$ is
also fibrewise linear, hence
\bea
\pa_u X(\eps,u)
&=&
\dd\wt{\Upsilon}_{(H(\eps,u)\cdot H(0,u)\inverse,\gamma_-(u))}
(\delta H_+,\Upsilon(a(0,u))_{\gamma_-(u)})\\
&+&
\dd\wt{\Upsilon}_{((H(\eps,u)\cdot H(0,u),\gamma_-(u))}
(\delta H_-,0_{\gamma_-(u)});
\eea
the first term of last expression vanishes, since it can be rewritten
as 
\bea\!\!\!\!\!\ \!\!\!\!\!
& &\!\!\!\!\!
\dd\wt{\Upsilon}_{@}
(\dd l_{(H(\eps,u)\cdot H(0,u)\inverse)}\dd\iota_\star a(0,u),
\Upsilon(a(0,u))_{\gamma_-(u)}\\
&=&\!\!\!\!\!
\dd\wt{\Upsilon}_{@}
(\dd l_{H(\eps,u)\cdot H(0,u)\inverse}\dd\iota_\star a(0,u),
\dd\wt{\Upsilon}_{(e_\star ,\gamma_-(u))}
(a(0,u),0_{\gamma_-(u)}))\\
&=&\!\!\!\!\!
\dd\wt{\Upsilon}_{@}
(\dd\mu^\star_{(H(\eps,u)\cdot H(0,u)\inverse,e_\star)}
(\dd l_{(H(\eps,u)\cdot H(0,u)\inverse)}\dd\iota_\star a(0,u),
a(0,u)),0_{\gamma_-(u)})\\
&=&\!\!\!\!\!
\dd\wt{\Upsilon}_{@}
(0_{H(\eps,u)\cdot H(0,u)\inverse},0_{\gamma_-(u)})=0
\eea
by equivariance, where we have set  
$@=(H(\eps,u)\cdot H(0,u)\inverse,\gamma_-(u))$  
to simplify the expressions; therefore
\bea
\pa_u X(\eps,u)
&=&
\dd\wt{\Upsilon}_{((H(\eps,u)\cdot H(0,u),\gamma_-(u))}
(\dd r_{H(\eps,u)\cdot H(0,u)\inverse}a(\eps,u),0_{\gamma_-(u)})\\
&=&
\dd\wt{\Upsilon}_{(e_\star,H(\eps,u)\cdot H(0,u)\inverse *\gamma_-(u))}
(a(\eps,u),0_{H(\eps,u)\cdot H(0,u)\inverse *\gamma_-(u)})\\
&=&\Upsilon(a(\eps,u))_{X(\eps,u)}\hs{3}.
\eea
We just have shown that
\be
\dd X_{(\eps,u)}
&=&
\Upsilon(\delta_r^\eps H(\eps,u))_{X(\eps,u)}\cdot d\eps 
+
\Upsilon(\delta_r^u H(\eps,u))_{X(\eps,u)}\cdot du\\ 
&=&
\Upsilon(h(\eps,u))_{X(\eps,u)}
\ee
and this concludes the proof.
\section{Notations and conventions} 
A bullet, ``$\bullet$'', usually denotes the one point manifold and
$\gr{f}=\{x,f(x)\}$ is the graph of a map $f$.
For any vector bundle $E$ and section $B\in\Gamma(\wedge^2 E)$
$B^\sharp:E^*\to E$ is the associated bundle map.
The structural map of a typical groupoid are source $\sor$, target
$\tar$, unit section $\eps$, inversion $\iota$ and partial
multiplication $\mu$; in special cases we use additional symbols,
which should be clear from the context.
We denote with   $\calG^{(n)}=\{(g_1\,,\,\dots\,,\,
g_n)\:|\:\sor(g_i)=\tar(g_{i+1})\}$, $n\in\NN$, the nerves of a
groupoid or a differentiable graph $(\calG,M;\sor,\tar)$, i.e. the 
strings of composable  elements.
As it is customary, we call a groupoid source 0- or  1-connected, if
it has 0- or 1-connected source fibres, respectively. 
The conormal bundle of a submanifold $C\subset M$ is
$N^*C\equiv\sf{Ann}\,TC\subset T^*M$; if $M$ is symplectic $T^\ohm C$
denotes the symplectic orthogonal bundle. 
Differently from \cite{st}, in this paper we denote with $\ol{P}$ the
opposite Poisson manifold $(P,-\pi)$ of a given Poisson manifold
$(P,\pi)$. 
The Hamilltonian vector field of $f\in\cif(M)$ is  $X^f=\{f,\cdot\}$
if $M$ is symplectic or Poisson; in order to make this  consistent we
define the Poisson bracket on $\cif(M)$ as $\{f,g\}=\ohm(X^g, X^f)$,
in the first case and $\{f,g\}=\pi(\dd f, \dd g)$ in the second.
\section*{Acknowledgments} 
I am deeply indebted to Kirill Mackenzie, for enlightening
conversations over the material presented in this paper and for his
precious comments during the writing stage, and with Alberto Cattaneo
for his constant support and encouragement. I also wish to thank Ping
Xu for stimulating discussions, Marco Zambon for his many valuable
suggestions on various versions of  this paper and Rui Loya Fernandes
for comments on an early version of it. 
%
%
%
%
%
%
%
%
%
%
%
%


\end{document}